\documentstyle{article}
\newcommand{\R}{I\!\! R}
\newcommand{\C}{I\!\! C}

\begin{document}

Tsemo Aristide

 ICTP,

Strada Costiera, 11

Trieste, Italy

tsemo@ictp.trieste.it

\bigskip

\bigskip

\bigskip

\bigskip

\centerline{\bf STRUCTURES GEOMETRIQUES.}

\bigskip

\bigskip

\bigskip

{\bf 0. Introduction.}

\medskip

Soient $X$ une vari\'et\'e diff\'erentiable connexe, $G$ un groupe
de Lie qui agit transitivement sur $X$ et dont l'action v\'erifie
le principe du prolongement analytique, ceci veut dire que deux
\'el\'ements de $G$ qui coincident sur un ouvert  de $X$
coincident sur $X$.
 Une structure g\'eom\'etrique model\'ee sur $X$, ou encore $(X,G)$
structure, est la donn\'ee d'une vari\'et\'e diff\'erentiable $M$,
d'un atlas $(U_i,\phi_i)$, o\`u $\phi_i:U_i\rightarrow X$ est un
diff\'eomorphisme
 sur son image tels que
$$
\phi_i\circ{\phi_j^{-1}}_{\phi_j(U_i\cap U_j)}:\phi_j(U_i\cap
U_j)\longrightarrow \phi_i(U_i\cap U_j)
$$
est la restriction d'un \'el\'ement $g_{ij}$ de $G$.

La famille $g_{ij}$ v\'erifie la relation de Chasles
$$
g_{ij}g_{jk}=g_{ik}.
$$
La $(X,G)$ structure  de $M$ se rel\`eve sur son rev\^etement
universel $\hat M$. Cette structure relev\'ee est d\'efinie par un
diff\'eomorphisme local $D:\hat M\rightarrow X$, et donne lieu \`a
une repr\'esentation du groupe fondamental de $M$:
$$
h_M:\pi_1(M)\longrightarrow G
$$
appel\'ee la repr\'esentation d'holonomie de la $(X,G)$ structure.

La d\'eveloppante peut aussi \^etre construite de la mani\`ere
suivante: on considere sur $M$, le faisceau des $(X,G)$
transformations locales \`a valeurs dans $X$, on note $F$ son
espace \'etal\'e, on a une application:
$$
F\longrightarrow X
$$
$$
[f]_x\rightarrow f(x)
$$
ou $x$ est un \'el\'ement de $M$, et $[f]_x$ un \'el\'ement de la
fibre de $x$, cette application est bien d\'efinie car $G$
v\'erifie le principe du prolongement analytique, et sa
restriction a une composante connexe de $F$ est la d\'eveloppante
(a un rev\^etement pr\`es).

\medskip

Des exemples de telles structures sont:

 les vari\'et\'es affines de
dimension $n$ o\`u $X$ est ${\R}^n$ et $G$ le groupe des
transformations affines de ${\R}^n$,

les vari\'et\'es projectives de dimension $n$ o\`u $X$ est
l'espace projectif $P{\R}^n$, et $G$ le groupe des transformations
projectives $PGL({\R}^n)$, ect...

Pour plus d'informations \`a propos de ces structures voir (3).

\medskip

On remarque que les ingr\'edients n\'ecessaires pour d\'evelopper
une th\'eorie des $(X,G)$ vari\'et\'es sont des invariants de
$1-$homotopie.

La $1-$homotopie des topos \'etant bien connue, on se propose de
 g\'en\'eraliser cette th\'eorie aux topos afin de traiter sur un
 m\^eme pied des applications en g\'eom\'etrie alg\'ebrique.

\medskip

\bigskip

{\bf 1. Rappels sur le groupe fondamental d'un topos.}

\medskip

La th\'eorie du groupe fondamental d'un topos est bien connue,
elle est propos\'ee en exercice dans 4. p. 321 Exercice 2.7.5.

\medskip

{\bf Definition 1.1.}

 - Soit $C$ une cat\'egorie, un crible de $C$ est une partie $R$
des objets $Ob(C)$ de $C$ telle que pour toute fl\`eche
$m:X\rightarrow Y$, si $Y\in R$, alors $X$ appartient \`a $R$.

Soient $f:C'\rightarrow C$ un foncteur et $R$ un crible de $C$, on note
$R^f$ le crible image r\'eciproque de $C$ c'est le crible de $C'$
 constitu\'e des objets dont l'image par $f$ appartient a $R$.

On note $C_X$ la sous-cat\'egorie de $C$ des objets au-dessus de
$X$.

\medskip

{\bf Definition 1.2.}

Une topologie de Grothendieck sur $C$ est une application qui a
tout objet $S$ associe une partie non vide $J(S)$ de l'ensemble
des cribles de la cat\'egorie $E_S$ au-dessus de $S$, telle que

pour toute fl\`eche $f:T\rightarrow S$ de $C$ et tout crible $R$
appartenant \`a $J(S)$, $R^f$ appartient \`a $J(T)$.

Pour tout objet $S$ de $C$,  tout \'el\'ement $R$ de $J(S)$, et
tout
 crible $R'$ de $E_S$, $R'$ appartient \`a $J(S)$ d\`es que pour
tout objet $f:T\rightarrow S$ de $R$, on a $R'^f$ appartient \`a
$J(T)$.

Les \'el\'ements de $J(S)$ sont appel\'es les raffinements de $S$.

Une cat\'egorie munie d'une topologie est appel\'ee un site.

\medskip

{\bf D\'efinitions 1.3.}

- Un pr\'efaisceau (d'ensemble) sur $C$ est un foncteur
contravariant d\'efinit sur $C$ et \`a valeurs dans la cat\'egorie
des ensembles.

- Un faisceau (d'ensemble) sur $C$ est un foncteur contravariant
$F$ de $C$ dans la cat\'egorie des ensembles tel que pour tout
raffinement $R$ de $S$, l'application
$$
F(S)\longrightarrow {\buildrel {lim F_R}\over{\longleftarrow}}
$$
soit bijective, o\`u $F_{\mid R}$ est le pr\'efaisceau d\'efini sur $R$
 par $F_{\mid R}(f)= F(T)$ pour toute fl\`eche
$f:T\rightarrow S$ de $R$.

\medskip

{\bf Exemple.}

On considere la cat\'egorie $e$ dont l'ensemble des objets est le
singleton  ayant pour seul \'el\'ement $x$, et l'ensemble des
morphismes de $x$ est r\'eduit a l'identit\'e.  Un faisceau de $e$
est la donn\'ee d'un ensemble.

\medskip

Soit $C$ une cat\'egorie, dont l'ensemble des objets $Ob(C)$ n'est
pas vide, il existe un foncteur projection $e_C:C\rightarrow e$
qui a tout objet associe $x$ et a toute fl\`eche l'identit\'e de
$x$.

\medskip

{\bf D\'efinitions 1.4.}

-Un faisceau constant $F$ de $C$ est un faisceau qui se factorise
par $e_C$. Autrement dit, il existe un faisceau $F_e$ de $e$ tel
que $F=F_e\circ e_C$.

- On dira que le faisceau $F$  de $C$ est localement constant si
et seulement si il existe une famille couvrante $(X_i)_{i\in I}$
de $C$ telle que la restriction de $F$ \`a la cat\'egorie
au-dessus de $X_i$, $C_{X_i}$ est constante.

\medskip

Dans la suite, on supposera que $C$ est un topos, ceci signifie
que $C$ est \'equivalente \`a la cat\'egorie des faisceaux
d\'efinit sur un  $U-$site standard (o\`u $U$ est un univers
donn\'e), ou de mani\`ere \'equivalente que les propri\'et\'es
suivantes sont v\'erifi\'ees:

(i) Munie de sa topologie canonique $C$ devient un $U-$site tel
que tout $U-$faisceau soit repr\'esentable.

(ii) Munie de sa topologie canonique, $C$ devient un $U-$site et
de plus

- Les limites projectives finies existent dans $C$,

- Les sommes index\'ees par un \'el\'ement de $U$ existent dans
$C$ et sont disjointes et universelles,

- Les relations d'\'equivalences sont \'effectives et
universelles.

\medskip

{\bf D\'efinition 1.5.}

Un morphisme de topos $f:X\rightarrow Y$, est un foncteur
$f^{-1}:Y\rightarrow X$ tel que:

- Pour tout faisceau sur $X$, le pr\'efaisceau
$f_*(F)(Y)=F(f^{-1}(Y))$ est un faisceau

- Le foncteur adjoint \`a gauche $f^*$ de $f_*$ commute aux
limites projectives finies.

\medskip

{\bf D\'efinition 1.6.}

Soit $Z$ l'objet final d'un topos $C$. On dira que $C$ (dont
l'ensemble des objets est non vide) est connexe, s'il n'existe pas
de famille couvrante $S$, constitu\'ee de deux \'el\'ements $S_1$
et $S_2$ telle que $S_1\times_Z S_2$ repr\'esente l'objet vide
(l'objet vide est l'objet qui repr\'esente le faisceau qui a un
objet associe l'ensemble vide).

\medskip

{\bf D\'efinition 1.7.}

- Une famille topologiquement couvrante $(X_i)_{i\in I}$ d'un
topos $C$ est connexe si pour tout $i$, $X_i$ est connexe.

\medskip

- Un topos  est localement connexe si et seulement si pour toute
famille topologiquement  couvrante $(X_i)_{i\in I}$, il existe une
famille
 topologiquement couvrante connexe $(Y_j)_{j\in J}$, telle que pour tout
 $j$, il existe un \'el\'ement $k(j)$ de $I$ tel que $Y_j$ est un sous-objet
 de $X_{k(j)}$.

 \medskip

{\bf D\'efinition 1.8.}

Soit $(X_i)_{i\in I}$ une famille topologiquement couvrante
connexe d'un topos. On appellera chemin de cette famille qu'on
note $(i_1,...,i_n)$, une famille d'\'el\'ements $(X_1,...,X_n)$
telle que $X_i\times_ZX_j$ est distinct de l'objet vide, ou $Z$
repr\'esente l'objet final.

\medskip

{\bf Proposition 1.9.}

{\it Un topos localement connexe est connexe si et seulement si, il existe une
famille topologiquement couvrante connexe $(X_i\rightarrow X)_{i\in I}$
telle que pour tout objet $X$ de $C$, pour tout $i$ de $I$, il existe
un chemin $(i_1=i,...,i_n)$ tel que $X_n\times_Z X$ est distinct de l'objet
vide ou $Z$ est l'objet final. }

\medskip

{\bf Preuve.}

Supposons qu'il existe un $X_i$ et un objet $X$ tel que pour tout
chemin issu de $X_i$, $X_n\times_ZX$ est l'objet vide. Soit
$C({X_i})$ l'ensemble des objets de $C$ tels que pour chacun de
ses \'el\'ements $Y$, il existe un chemin $(i_1=i,...,i_n)$ tel
que $X_n\times_ZY$ soit distinct de l'objet vide. Son
compl\'ementaire $U_i$ est non vide, l'objet final est alors la
somme directe de la somme de l'ensemble des objets de $U_i$ de la
somme de l'ensemble des objets de  $C(X_i)$, d'o\`u le r\'esultat.

\medskip

 On appelle $C({X_i})$ la composante connexe de $X_i$. Un topos localement
 connexe est somme directe de ses composantes connexes.

Dans la suite on supposera que les topos consid\'er\'es sont
localement connexes et connexes.

\medskip

On munit l'ensemble des chemins associ\'es \`a la famille
topologiquement couvrante connexe $(X_i)_{i\in I}$ de la relation
d'\'equivalence suivante:

Deux chemins $x$ et $y$ sont \'equivalents si et seulement si il
existe une suite de chemins $z_1,...,z_n$ tels que pour $k\leq
n-1$, $z_k=(i_1,...,i_l,i_l,i_{l+1},..,i_{m})$ et
$z_{k+1}=(i_1,..,i_l,i_{l+1},..i_m)$ avec $x=z_1$, $y=z_n$ ou
$x=z_n$, $y=z_1$.

On note $Chem((X_i)_{i\in I})$ l'ensemble des classes
d'\'equivalences.

Dans la classe d'\'equivalence de tout chemin $x$, il existe un
chemin $\bar x=(j_1,...,j_k)$ tel que $j_r$ est distinct de
$j_{r+1}$.

\medskip

Soient $x=(i_1,...,i_n)$ et $y=(j_1,...,j_m)$ deux chemins
repr\'esentant les \'el\'ements respectifs $\bar x$ et $\bar y$ de
 $Chem((X_i)_{i\in I})$ tels que $i_n=j_1$, on associe a $\bar x$ et $\bar y$
 l'\'el\'ement $\bar{x*y}$ de $Chem((X_i)_{i\in I})$ dont un repr\'esentant est
 $x*y=(x_1,..,x_n,j_2,..j_n)$.

 On peut maintenant d\'efinir le groupoide  $Gr((X_i)_{i\in I})$
 dont les objets sont les $i$, o\`u
 $X_i$ est un \'el\'ement de la famille couvrante $(X_i)_{i\in I}$.

 Un morphisme entre $i$ et $j$ est la classe d'un chemin $(i_1,...,i_n)$ telle
 que $i_1=i$ et $i_n=j$.

L'inverse de la classe de $(i_1,...,i_n)$ est $(i_n,i_{n-1},..,i_1)$.

L'ensemble des morphismes repr\'esent\'es par des \'el\'ements de la forme
$(i,i_2,..,i_{n-1},i)$ est un groupe qu'on note $Aut(i)$.

 \bigskip

Soit $F$ un faisceau localement constant d\'efini sur $C$,
$(X_i)_{i\in I}$ une famille topologiquement couvrante connexe
 de $C$ telle que la restriction  de $F$ a chaque
$X_i$ est un faisceau constant. Une telle famille sera appel\'ee famille
trivialisante.

\medskip

Soient $X$ et $Y$ deux objets de $C$ tels que la restriction de
$F$ \`a $C_X$ et $C_Y$ est constante et $Z$ l'objet final de $C$.
Supposons que $X\times_ZY$ n'est pas l'objet initial, alors les
projections $p_x:X\times_ZY\rightarrow X$ et
$p_y:X\times_ZY\rightarrow Y$ donnent lieu \`a deux isomorphismes
d'ensembles $g_x:F(X)\rightarrow F(X\times_ZY)$ et
$g_y:F(Y)\rightarrow F(X\times_ZY)$. On en d\'eduit un
isomorphisme $g_{xy}=g_x{g_y}^{-1}:g_y(F(Y))\rightarrow
g_x(F(X))$. On choisit pour $g_{xx}$ l'application identit\'e.

Appliquons ceci au chemin $(i_1,...,i_n)$, on en d\'eduit pour
$i_k$ un isomorphisme $g_{i_ki_{k+1}}:
g_{i_{k+1}}(F(X_{i_{k+1}}))=F(X_{i_k}\times_Z
X_{i_{k+1}})\rightarrow g_{i_k}(F(X_{i_k}))$, et par suite un
isomorphisme $g_{i_1i_n}:g_{i_n}(F(X_{i_n}))\rightarrow
g_{i_1}(F(X_{i_1}))$. $g_{i_1i_n}= g_{i_1i_2}\circ..\circ
g_{i_{n-1}i_n}$.

 On en d\'eduit une repr\'esentation $hol_F:Aut(i)\rightarrow Aut(F(X_i))$ ou
 $Aut(F(X_i))$ d\'esigne le groupe des automorphismes de l'ensemble
 $F(X_i)$.
 On pose $F_S(X_i)=hol_F(Aut(i))$ son image.

 \medskip

\medskip

 {\bf Proposition 1.10.}

{\it Supposons que le topos localement connexe $C$ soit connexe,
alors pour tous objets $X_i$, $X_j$ de la famille topologiquement
couvrante connexe $(X_i)_{i\in I}$, les groupes $hol_F(Aut(i))$ et
$hol_F(Aut(j))$ sont isomorphes et de plus ne d\'ependent pas de
la famille couvrante trivialisante connexe choisie.}

\medskip

{\bf Preuve.}

Soient $X_i$ et $X_j$ deux \'el\'ements de $(X_i)_{i\in I}$, il
existe un chemin $(i_1,...,i_n)$ entre $X_i$ et $X_j$, la classe
de ce chemin induit un isomorphisme entre $hol(Aut(i))$ et
$hol(Aut(j))$.

Montrons maintenant que la classe d'isomorphisme de $hol(Aut(i))$ ne
d\'epend pas de la famille topologiquement couvrante connexe trivialisante
de $F$
 choisie.

 Consid\'erons une autre famille topologiquement couvrante connexe
 trivialisante $(Y_j)_{j\in J}$ de $F$. Il existe une famille
 topologiquement couvrante connexe trivialisante $(Z_k)_{k\in K}$ telle
 que  tout $Z_k$, est un sous-objet d'un objet $X_{i(k)}$
 de $(X_i)_{i\in I}$, et d'un objet $Y_{j(k)}$ de $(Y_j)_{j\in J}.$

 Il suffit de montrer que $hol_F(Aut(X_i))$ et $hol_F(Aut(Z_k))$ sont isomorphes,
 pour toute famille trivialisante $(Z_k)_{k\in K}$ telle qu'il existe
 $i(k)$ tel que $Z_k$ soit un sous-objet de $X_{i(k)}$ pour tout $k$.

 Soit $x=(k_1=k,...,k_n=k)$ un chemin de $(Z_k)_{k\in K}$, On consid\`ere
 $X_{i_l}$ contenant $Z_{k_l}$, $y=(i_1,...,i_n)$ est un chemin de
 $(X_i)_{i\in I}$ tel que $X_{i_1}=X_{i_n}$.  $hol_F(x)=hol_F(y)$, on d\'efinit ainsi un morphisme de
 groupe de $hol_F(Aut(Z_k)_{k\in K})$ dans $hol_F(Aut(X_i)_{k\in I}).$

 D\'efinissons son inverse, consid\'erons un chemin $(i_1=i,...,i_n=i)$ de
 $(X_i)_{i\in I}$. Pour chaque $X_{i_l}$, il existe un de ses
  sous-objets $Z_{k_l}$
 appartenant \`a la famille $(Z_k)_{k\in K}$ tel que $Z_{k_l}\times_ZX_{k_{l+1}}$ soit diff\'erent de
 l'objet vide. On peut relier $Z_{k_l}\times_ZX_{k_{l+1}}$
 et $Z_{k_{l+1}}$ par un chemin de $X_{k_{l+1}}$ car il est
 connexe. On obtient ainsi un chemin de $(Z_k)_{k\in K}$ qui
 d\'efinit l'inverse du morphisme pr\'ec\'edent.

 \bigskip

 On va noter $(G_i)_{i\in I}$ l'ensemble des groupes $G_i$ tels que
 $G_i=hol_F(Aut(X_i))$, o\`u $X_i$ est un \'el\'ement d'une famille
 topologiquement couvrante trivialisante du faisceau $F$.

 On dira que $G_i\leq G_j$ s'il existe un morphisme surjectif entre
 $G_j$ et $G_j$, on obtient ainsi un systeme projectif dont le
 compl\'et\'e projectif est appel\'e progroupe fondamental du topos connexe
 et localement connexe $C$. On le note $pro\pi_1(C)$.

 \medskip

 Pour un topos localement connexe, on peut d\'efinir le progroupe fondamental
 de chacune de ses composantes connexe.

 \medskip

 {\bf D\'efinition 1.11.}

 On dira qu'un topos  connexe $C$
 est simplement connexe si et seulement si tout faisceau localement
 constant d\'efini sur $C$ est
 constant.

- On dit que $C$ est localement simplement connexe, si et
seulement si de toute famille topologiquement couvrante, on peut
extraire une famille topologiquement couvrante  $(X_i)_{i\in I}$
telle que $X_i$ est connexe et simplement connexe.

\medskip

{\bf Remarque.}

\medskip

 Pour un topos connexe et localement simplement connexe, le progroupe
 fondamental est un groupe qui se d\'ecrit comme suit:

Soit $(X_i)_{i\in I}$ une famille topologiquement couvrante
connexe et simplement connexe,
 soit $G$ l'intersection des noyaux des
 repr\'esentations $hol_F:Aut(i)\rightarrow Aut(F(X_i))$ pour $i$ fix\'e et $F$
 un faisceau localement constant de $C$ quelconque.
 le groupe fondamental est isomorphe au quotient de $Aut(i)$ par $G$.

 \bigskip

{\bf 2. Structures g\'eom\'etriques en th\'eorie des topos.}

 \medskip

Le but de cette partie est d'\'etendre la notion de structure
g\'eom\'etrique \`a la th\'eorie des topos.

\medskip

 {\bf D\'efinition 2.1.}

Soit $C$ un topos, $G$ un sous-groupe d'automorphismes de $C$,
on dira que l'action de $G$ v\'erifie le principe du prolongement
analytique (p.p.a) si et seulement si deux \'el\'ements de $G$ qui coincident
sur la cat\'egorie au-dessus d'un objet de $C$ coincident sur $C$.

\medskip

Dans la suite les topos consid\'er\'es seront   localement
connexes, et on supposera qu'ils  ont un nombre fini de
composantes connexes.

\medskip

{\bf D\'efinition 2.1.}

Soit $C$ un topos, $G$ un groupe d'automorphismes de $C$ qui
v\'erifie le p.p.a. On dira qu'un topos $D$ est muni d'une $(C,G)$
structure si et seulement si

- il existe une famille topologiquement couvrante connexe de $D$,
$(X_i)_{i\in I}$, telle que pour tout $X_i$, il existe un
isomorphisme local
 de topos $\phi_i:X_i\rightarrow C$
 i.e il existe un sous-objet $Y_i$ de $C$ tel que $\phi_i$ se
 factorise par un isomorphisme
 $$
 \phi_i:X_i\rightarrow Y_i.
 $$

- Supposons que $X_i\times_ZX_j$ soit distinct de l'objet initial
($Z$ est l'objet final), il existe un \'el\'ement  $g_{ij}$ de $G$
tel que
$$
g_{ij}(\phi_j)_{X_i\times_ZX_j}= (\phi_i)_{X_i\times_ZX_j}
$$
ou $(\phi_i)_{X_i\times_ZX_j}$ (resp $(\phi_j)_{X_i\times_ZX_j}$)
est la restriction de $\phi_i$ \`a la cat\'egorie au-dessus de
$X_i\times_ZX_j$, (resp. $\phi_j$ \`a la cat\'egorie au-dessus de
$X_i\times_ZX_j$).

\medskip

On a
$$
g_{jk}{(\phi_k)}_{X_i\times_ZX_j\times_Z X_k}=
(\phi_j)_{X_i\times_ZX_j\times_Z X_k},$$

$$
g_{ij}(\phi_j)_{X_i\times_ZX_j\times_Z X_k}=
{(\phi_i)}_{X_i\times_ZX_j\times_Z X_k}.$$

On en d\'eduit que

$$
g_{ij}g_{jk}(\phi_k)_{X_i\times_ZX_j\times_Z X_k}=
({\phi_i)}_{X_i\times_ZX_j\times_Z X_k}$$ et par suite que
$$
g_{ij}g_{jk}=g_{ik} \leqno (1)
$$
d'apr\`es le principe du prolongement analytique.

\medskip

La relation $(1)$ permet de d\'efinir un faisceau localement
constant $F$ sur $D$, tel que  la restriction de $F$ \`a la
cat\'egorie au-dessus de  $X_i$ est $G$, et pour tout $X$, $F(X)$
est le noyau du couple de fl\`eches
$$
\prod F(X_i){\buildrel{\rightarrow}\over{\rightarrow}}\prod
F(X_i\times_XX_j),
$$
Les morphismes de transitions de $F$ sont induits par les
$g_{ij}$.

On en d\'eduit une repr\'esentation du progroupe fondamental
$pro\pi_1(D_0)$ d'une composante connexe $D_0$ de $D$ dans $G$,
donn\'e par la projection

$$
hol_{D_0}:pro\pi_1(D_0)\longrightarrow hol_F((X_i)_{i\in I})
$$
appel\'ee repr\'esentation d'holonomie de la $(C,G)$ structure de $D_0$.

\medskip

{\bf D\'efinition 2.2.}

Soit $D$ et $D'$ deux $(C,G)$ topos, on dira qu'un isomorphisme local
$f:D\rightarrow D'$ est un $(C,G)$ morphisme si et seulement si

$$
\phi_{i(j)}f_{\mid X_i}=k_{ij}\phi_i
$$

o\`u $k_{ij}$ est un \'el\'ement de $G$, et $(X_i)_{i\in I}$,
$(Y_j)_{j\in J}$ des familles couvrantes d\'efinissant les $(C,G)$
structures de $D$ et $D'$ et $f_{\mid X_i}$ la restriction de $f$
\`a la cat\'egorie au-dessus de $X_i$. On suppose sans restreindre
la g\'en\'eralit\'e que $f_{X_i}$ se factorise par un sous-objet
de $Y_{i(j)}$. On munit ainsi l'ensemble des $(C,G)$ topos d'une
structure de cat\'egorie.

\medskip

L'existence d'un rev\^etement universel n'\'etant pas sure, on le
remplace par le faisceau $H_S$ des $(C,G)$ morphismes locaux de
$D$.

A $H_S$, on associe la cat\'egorie dont les objets sont les
couples $(X,s)$, ou $X$ est un objet de $D$ et $s$ un \'el\'ement
de $F(X)$, une fl\`eche entre deux objets $(X,s)$ et $(Y,s')$ de
$H$, est une fl\`eche $f:X\rightarrow Y$ telle que $H_S(s')=s$.
Cette cat\'egorie est un topos.

On d\'efinit alors le morphisme:
$$
Dev: H\rightarrow C
$$
d\'efini localement  sur la cat\'egorie au-dessus de $(X,s)$ par
$s$

Le morphisme $s$ est un morphisme de topos entre $D_X$ et $C$, il
induit un morphisme de topos sur la cat\'egorie au dessus de
$(X,s)$.

L'application $Dev$ est appel\'ee d\'eveloppante de la $(C,G)$
structure de $D$.

\bigskip

On consid\`ere maintenant $(C,G)$  (resp.$(C',G')$) un topos $C$
(resp. $C'$) et $G$ (resp. G') un groupe de $C$ (resp. $C'$) qui
v\'erifie le p.p.a. On suppose qu'on a un isomorphisme local de
topos $\phi: C\rightarrow C'$ et une repr\'esentation
$\Phi:G\rightarrow G'$ telle que pour tout \'el\'ement $g$ de $G$,
$\phi \circ g=\Phi(g)\circ\phi$.

A tout $(C,G)$ topos $D$, d\'efinit par la famille couvrante
$((X_i)_{i\in I},\phi_i)$, on peut associer le $(C',G')$ topos d\'efini par
 $((X_i)_{i\in I},\phi\circ \phi_i)$.

 On d\'efinit ainsi un foncteur $\phi_*$ de la cat\'egorie des $(C,G)$ topos
 dans la cat\'egorie des $(C',G')$ topos.
\bigskip

{\bf 3. Espace des $(C,G)$ structures.}

\medskip

On consid\`ere un topos $D$ muni d'une $(C,G)$ structure, soit
$((X_i)_{i\in I},\phi_i)$ la famille trivialisante qui la
d\'efinit

Consid\'erons le topos $S((X_i)_{i\in I})$ dont l'objet final
 est obtenu en recollant les produits de topos
$U_i\times C$ par la relation de Chasles
$$
\tilde g_{ij}:U_i\times_ZU_j\times C\rightarrow U_i\times_ZU_j\times C
$$
induit par les fonctions de transitions  $g_{ij}$.

$S((X_i)_{i\in I})$ est appel\'e le topos structural de la
 $(C,G)$ structure de $D$.

\medskip

Les projections $p_i:U_i\times C\rightarrow U_i$ se recollent en une
projection $p:S\rightarrow D$ car les fonctions de transitions
$\tilde g_{ij}$ induisent l'identit\'e sur le premier facteur.

\medskip

{\bf D\'efinition 3.1.}

Un morphisme de topos $s:D\rightarrow S((X_i)_{i\in I})$ telle que
$p\circ s= Id_D$ sera appel\'e une section du topos
$S((X_i)_{i\in I})$.

Une section est transverse si et seulement si pour tout $i$,
$p_i\circ s_i$ est un isomorphisme  local. ($s_i$ est la
restriction de $s$ a $U_i$).

Pour une $(C,G)$ structure $D$, on fixera une section transverse
$s_0$ de $S((X_i)_{i\in I}$.

\medskip

{\bf Proposition 3.2.}

{\it La donn\'ee d'une section transverse permet de d\'efinir une $(C,G)$
 structure sur $D$ par $((X_i)_{i\in I},p_is_i)$.}

\medskip

L'espace des $(C,G)$ structures de $D$ peut aussi \^etre d\'ecrit comme
l'ensemble des triples suivants:

- La donn\'ee d'une $(C,G)$ structure  de  $D$, d\'efinit par la
famille trivialisante $((X_i)_{i\in I},\phi_i)$, on fixe une carte
$(X_{i_0},\phi_{i_0})$,

- d'un isomorphisme de
topos $D'\rightarrow D$.

\medskip

On note $S(D,C,G)$ l'ensemble des $(C,G)$ structures de $D$.

On a une application
$$
S(D,C,G)\longrightarrow Hom(pro\pi_1(D_0),G)
$$
qui a un \'el\'ement de $S(D,C,G)$ associe sa repr\'esentation d'holonomie.

\bigskip

\bigskip

{\bf 4. Applications \`a la g\'eom\'etrie alg\'ebrique.}

\medskip

On consid\`ere un sch\'ema $S$  d\'efinit sur un corps $k$, un
sous-groupe $G$ du groupe des automorphismes
 de $S$ v\'erifiant le p.p.a,
 (on peut prendre $S$ irr\'eductible et $G$ un sous-groupe du groupe
 des automorphismes
 de $G$).
  Il op\`ere aussi
sur le site \'etale $Sit_S$ de $S$. On peut d\'efinir donc une
notion de $(S,G)$ sch\'emas dans la cat\'egorie des sch\'emas
d\'efinis sur $k$.

Un $(S,G)$ sch\'ema est un sch\'emas $T$, dont le site \'etale
$Sit_T$ est muni d'une structure de $(S,G)$ structure.

Supposons  que le groupe $hol(T)$ d'holonomie de la $(S,G)$ structure soit finie,
 alors il d\'efinit un faisceau en groupe finis au-dessus de $T$ qu'on note $F$.
 On sait d'apr\`es la th\'eorie du groupe fondamental de Grothendieck qu'il
 existe un sch\'ema $\hat T$ \'etale au-dessus de $T$
 tel que le relev\'e de $F$ sur $T$ soit trivial.
 on l'appelle un rev\^etement d'holonomie de la $(S,G)$ structure de $T$.

 Les ouverts \'etales qui d\'efinissent la $(S,G)$ structure de $T$ se
  rel\`event
 sur $\hat T$ (par l'op\'eration produit fibr\'e par $\hat T$)
 et d\'efinissent un mophisme de sch\'emas de $\hat T\rightarrow S$
 qui est la d\'eveloppante de cette $(S,G)$ structure.

La   repr\'esentation d'holonomie d\'efinit une repr\'esentation
$$
hol:pro\pi_1(T)\rightarrow hol(T)
$$
ou $pro\pi_1(T)$ est le progroupe fondamental de $T$ au sens de la
g\'eom\'etrie alg\'ebrique.

\medskip

Consid\'erons une $(S,G)$ structure d\'efinie sur le sch\'ema $T$
de groupe d'holonomie fini $h(T)$, on note $\hat T$ le
rev\^etement correspondant. Pour $\hat T$ fix\'e, l'espace des
d\'eformations est l'ensemble des espaces  fibr\'es $D(T,\hat
T,S,G)$ d\'efinis par
$$
\gamma:\hat T\times S\longrightarrow \hat T\times S
$$
$$
(x,y)\longrightarrow (x,hol(T)(\gamma(y))
$$
ou $hol(T)$ est l'image de la repr\'esentation d'holonomie d'une
$(S,G)$ structure et $\gamma$ un de ses \'el\'ements.

Pour un tel fibr\'e
 on a sur cet espace un feuilletage horizontal horizontal
${\cal F}$ qui est le projet\'e du feuilletage $\hat{\cal F}$ de
$\hat T\times S$ dont les feuilles sont les sous-sch\'emas $\hat
T\times {y}$, ou $y$ est un point de $S$.

La donn\'ee d'une section $s$ de $T$ transverse \`a ce feuilletage
horizontal d\'efinie la $(S,G)$ structure de $T$. En effet elle se
rel\`eve en une section de $\hat s:\hat T\rightarrow \hat T\times
S$ qui se projette sur $S$ en un isomorphisme  local \'equivariant
pour la repr\'esentation d'holonomie.

Si $T'$ est un autre rev\^etement fini de $T$ tel qu'il existe une
surjection $T'\rightarrow\hat T$, on a une injection  de $D(T,\hat
T, S,G)\rightarrow D(T,T',S,G)$. La limite inductive des espaces
$D(T,\hat T,S,G)$ est l'espace classifiant des $(S,G)$ structures
de $T$, on le note $D(T,S,G)$.

\medskip

Supposons que les sch\'emas $S$ et $T$ soient d\'efinis sur
${\C}$, on munit le groupe $G$ de la topologie compacte-ouverte
induite par la topologie d'espace analytique de $S$. Dans ce cadre
on peut parler de d\'eformations des repr\'esentations du groupe
fini $hol(T)$ dans $G$.

\medskip

{\bf Lemme 4.1.}

{\it Muni de la topologie qu'on vient de d\'efinir, Soit
$(\rho_t)_{t\in I}$ une d\'eformation continue d'une
repr\'esentation $\rho$ de $hol(T)$, ($I$ est un intervalle
contenant $0$ tel que $h_0=hol(T)$) alors, pour tout $t$, $\rho_t$
est isomorphe \`a $\rho$ o\`u $t$ appartient a un intervalle
ouvert contenu dans $I$ et contenant $0$. }

\medskip

{\bf Preuve.}

Soient $g_1,...,g_n$ les \'el\'ements de $hol(T)$. Consid\'erons
l'ensemble $X$ des \'el\'ements de $S$ tels que le cardinal de
$(g_1(x),...,g_n(x))=n$, $X$ est un ouvert non vide. Soit $x$ un
de ses \'el\'ements, il existe des ouverts $U_1$,...,$U_n$ (pour
la topologie induite par la structure analytique) contenant
respectivement $g_1(x)$,...,$g_n(x)$, deux \`a deux disjoints, des
ouverts $V_1$,...,$V_n$ contenant respectivement
$g_1(x)$,...,$g_n(x)$, et contenus respectivement dans
$U_1$,...,$U_n$, un intervalle $J$ contenant $0$ contenu dans $I$,
tels que pour tout $t$ appartenant a $J$, $\rho_t(g_i)(x)$ est un
\'el\'ement de $V_i$ et si $g_ig_k(x)$ a $V_l$, alors
$\rho_t(g_i)\rho_t(g_k)(x)$ appartient a $U_l$. On en d\'eduit que
l'application de $rho_t(hol(T))\rightarrow hol(T)$ qui a
$\rho_t(g_i)$ associe $g_i$ est un isomorphisme de groupe (t$\in
J)$.

\medskip

{\bf Th\'eor\`eme 4.2.}

{\it Supposons que $S$ et $T$ soint deux sch\'emas sur ${\C}$,
alors l'image de $D(S,T,\hat T,G)\rightarrow Hom(hol(T),G)$ est
ouverte.}

\medskip

{\bf Preuve.}

Sous les conditions du th\'eor\`eme, pour toute d\'eformation
$(\rho_t)_{t\in I}$ de la repr\'esentation d'holonomie de $T$, il
existe un intervalle $J$ inclu dans  $I$ tel que les fibr\'es
$D(T,\hat T,S,G)$ et le fibr\'e plat induit par les
repr\'esentations $\rho_t$, $t\in J$ sont isomorphes. On les
identifie, et on note ${\cal F}_{\rho_t}$ le feuilletage du
fibr\'e plat induit par $\rho_t, t\in J$ vu dans ce fibr\'e. Pour
tout $t$ suffisament petit, on peut supposer qu'il reste
transverse  a la section qui d\'efinit la $(S,G)$ structure
initiale, cette section  d\'efinit une $(S,G)$ structure
d'holonomie $\rho_t$.

\bigskip

{\bf Remarque.}

\medskip

Si on suppose que les morphismes de transitions $g_{ij}$ qui
permettent de d\'efinir la $(S,G)$ structure sont des submersions,
on peut d\'evelopper une telle th\'eorie itou.

\medskip

{\bf Bibliographie.}

\medskip

1. Deligne, P. Le groupe fondamental de la droite projective moins
trois points. Math. Sci. Res. Inst. Publ, 16.

2. Giraud, J. Cohomologie non ab\'elienne. Springer 1971.

3. Goldman, W. Geometric structures on affine manifolds and
varieties of representations. Contemo. Math, 74.

4. S.G.A.1 S\'eminaire dirig\'e par A.Grothendieck, Rev\^etement
\'etales et groupe fondamental. Springer 1971.

5. S.G.A.4 S\'eminaire dirig\'e par A. Grothendieck, Th\'eorie des topos
et cohomologie \'etales des sch\'emas. Spinger 1972

\end{document}